\documentclass[11pt]{amsart}
\usepackage{graphicx}
\usepackage{verbatim}
\usepackage{textcomp}
\usepackage{amssymb}
\usepackage{cite}
\usepackage{amsmath}
\usepackage{latexsym}
\usepackage{amscd}
\usepackage{amsthm}
\usepackage{mathrsfs}
\usepackage{xypic}
\usepackage{bm}
\usepackage{url}
\usepackage{hyperref}

\newtheorem{thm}{Theorem}%[section]

\newtheorem{lem}[thm]{Lemma}

\theoremstyle{definition}

\theoremstyle{remark}

\def\R{\mathbb R}

\def\H{\mathbb H}
\def\SS{\mathbb S}

\def\f{\frac}

\def\ra{\rightarrow}
\def\pt{\partial}

\begin{document}
\title[Non-collapsing for hypersurface flows]{Non-collapsing for hypersurface flows in the sphere and hyperbolic space}
\author{Ben Andrews}
\address{Mathematical Sciences Institute, Australia National University; Mathematical Sciences
Center, Tsinghua University; and Morningside Center for Mathematics, Chinese Academy of Sciences.}
\email{Ben.Andrews@anu.edu.au}
\author{Xiaoli Han}
\address{Department of mathematical sciences, Tsinghua University, 100084, Beijing, P. R. China}
\email{xlhan@math.tsinghua.edu.cn}
\author{Haizhong Li}
\address{Department of mathematical sciences, and Mathematical Sciences
Center, Tsinghua University, 100084, Beijing, P. R. China}
\email{hli@math.tsinghua.edu.cn}
\author{Yong Wei}
\address{Department of mathematical sciences, Tsinghua University, 100084, Beijing, P. R. China}
\email{wei-y09@mails.tsinghua.edu.cn}
%\date{}
\thanks{The research of the first author was partially supported by Discovery Projects grant DP120102462 of the Australian Research Council.
The research of the second and third author was supported by NSFC No. 11131007 and No. 10971110, respectively.}

\maketitle

\begin{abstract}
We prove a non-collapsing property for curvature flows of embedded hypersurfaces in the sphere and in hyperbolic space.
\end{abstract}

\section{Introduction}

Let $X: M^n\times [0,T)\ra (N^{n+1}(c),\bar{g})$ be a family of embedded hypersurfaces in the simply connected space-form $N^{n+1}(c)$ with sectional curvature $c$, evolving by the curvature flow
\begin{align}\label{F-flow}
    \f{\pt X(x,t)}{\pt t}=-F(x,t)\nu(x,t),
\end{align}
where $\nu$ is the unit outward normal, and the speed $F$ is given by a homogeneous degree one, monotone increasing function of the principal curvatures defined on a symmetric convex cone $\Gamma\subset\R^n$. We further assume that $F$ is normalized so that $F(1,\cdots,1)=n$.  The purpose of this paper is to prove a non-collapsing result for the flow \eqref{F-flow} in $N^{n+1}(c)$. Here $c=1$ corresponds to the sphere $\SS^{n+1}=\{X\in\R^{n+2}:\langle X,X\rangle=1\}$, and $c=-1$ corresponds to the hyperbolic space $\H^{n+1}$. We use the hyperboloid model of $\H^{n+1}$, i.e., $\H^{n+1}$ is the upper sheet $(x_0>0)$ of the two-sheeted hyperboloid $\{X:\ \langle X,X\rangle = -1\}$ in the Minkowski space $\R^{n+1,1}$.  In these expressions the inner product $\langle\cdot,\cdot\rangle$ refers to the inner product in $\R^{n+2}$ or $\R^{n+1,1}$ respectively.

Following \cite{ALM}, we define the function $k(x,y,t)$ for $y\neq x$ by
\begin{align}\label{defn-k}
    k(x,y,t)=&\f{2}{d^2}\langle X(x,t)-X(y,t),\nu(x,t)\rangle,
\end{align}
where $d=\|X(x,t)-X(y,t)\|$ is the distance. The supremum of $k(x,y,t)$ over $y$ gives the geodesic curvature of the largest interior sphere which touches at $x$. We call $\bar{k}(x,t)=\sup\{k(x,y,t):y\in M,y\neq x\}$ the interior sphere curvature at the point $(x,t)$, and $\underline{k}(x,t)=\inf\{k(x,y,t):y\in M,y\neq x\}$ the exterior sphere curvature at $(x,t)$. Note that the definitions of $\bar{k}(x,t)$ and $\underline{k}(x,t)$ involve extrema of $k(x,y,t)$ over the noncompact set $y\in M: y\neq x\}$.
The function $k$ extends continuously to a compactification obtained by adjoining the unit sphere in the tangent space at each point $(x,t)$ (see \cite{ALM}), and it follows that $\bar{k}(x,t)$ is no less than the maximum principal curvature $\kappa_{\max}(x,t)$, and either there exists $\bar{y}\in M\setminus\{x\}$ such that $\bar{k}(x,t)=k(x,\bar{y},t)$, or there exists a unit vector $v\in T_xM$ such that $\bar{k}(x,t)=h_{(x,t)}(v,v)=\kappa_{\max}(x,t)$. Similarly, $\underline{k}(x,t)$ is no greater than the minimum principal curvature $\kappa_{\min}(x,t)$.

In geometric flow problems, the idea of dealing with a function on the product $M\times M$ first appears in Huisken \cite{Hu98} and Hamilton's \cite{Ha2,Ha3} work on the curve shortening flow and Ricci flow. See also the recent refinements of these works by the first author and Paul Bryan \cite{AB1,AB2,AB3}. The first author \cite{An} used an argument of this kind to give a direct proof of Sheng-Wang's non-collapsing theorem \cite{SW} for mean-convex mean curvature flow in $\R^{n+1}$. Later this was generalized to fully nonlinear curvature flows by the first author, Langford and McCoy \cite{ALM}. Recently, the technique of \cite{An} was used by Brendle \cite{Bren} to prove the Lawson conjecture, and subsequently by the first and third authors \cite{AL} to prove the Pinkall-Sterling conjecture. In this paper, we follow the ideas in \cite{An,ALM} to prove the following non-collapsing properties of the flow \eqref{F-flow} in $\SS^{n+1}$ and $\H^{n+1}$.

\begin{thm}\label{sphere}
Let $X:M^n\times[0,T)\ra\SS^{n+1}$ be an embedded solution of \eqref{F-flow}. If $F$ is concave and positive, then we have $ \f{\bar{k}(x,t)}{F(x,t)}-\f 1{n}\leq C_1e^{-2nt}$ with $C_1=\sup\{\f{\bar{k}(x,0)}{F(x,0)}-\f 1{n}:x\in M\}\geq 0$. If $F$ is convex and positive, then $ \f{\underline{k}(x,t)}{F(x,t)}-\f 1{n}\geq C_2e^{-2nt}$ with $C_2=\inf\{\f{\underline{k}(x,0)}{F(x,0)}-\f 1{n}:x\in M\}\leq 0$.
\end{thm}

An important examples of flows of the form \eqref{F-flow} is the mean curvature flow, in which $F$ is equal to the mean curvature $H=\sum_i\kappa_i$ (here $\kappa_1,\dots,\kappa_n$ are the principal curvatures).  Since this is both concave and convex as a function of the principal curvatures, we get from Theorem \ref{sphere} that under the mean-convex embedded mean curvature flow in $\SS^{n+1}$, the following pinching result holds.
\begin{align*}
    C_2e^{-2nt}+\f 1n\leq \f{\underline{k}(x,t)}{H(x,t)}\leq \f{\bar{k}(x,t)}{H(x,t)}\leq \f 1n+C_1e^{-2nt}.
\end{align*}

In the hyperbolic case, we state a result only for the mean curvature flow:   We prove that

\begin{thm}\label{hbolic}
Let $X:M^n\times[0,T)\ra\H^{n+1}$ be an embedded solution of the mean curvature flow.
\begin{enumerate}
  \item\label{item1} If $M_0=X(M,0)$ is mean-convex, then $\f{\bar{k}(x,t)}{H(x,t)}-\f 1{n}\leq C_3e^{2nt}$ with $C_3=\sup\{\f{\bar{k}(x,0)}{H(x,0)}-\f 1{n}:x\in M\}$, and $\f{\underline{k}(x,t)}{H(x,t)}-\f 1{n}\geq C_4e^{2nt}$ with $C_4=\inf\{\f{\underline{k}(x,0)}{H(x,0)}-\f 1{n}: x\in M\}$.
  \item If $M_0=X(M,0)$ satisfies $H(x,0)>n$, then $\f{\bar{k}(x,t)}{H(x,t)-n}\leq C_5$ with $C_5=\sup\{\f{\bar{k}(x,0)}{H(x,0)-n}: x\in M\}$.
\end{enumerate}
\end{thm}

The result of case \ref{item1} allows the `collapsing ratio' to grow exponentially with time.  This should be expected.   For example, consider a non-compact convex region with two boundary components each having constant principal curvatures less than $1$:  The boundary sheets move together with exactly such a rate of non-collapsing.  Compactifying this example by intersecting with a large sphere and pasting in a large semi-cylindrical shell will produce compact examples for which the collapsing ratio decays in this way on arbitrarily long finite time intervals.

Despite the fact that the collapsing ratio can become large, the result still provides useful information concerning finite time singularities.
An inspection of the proof shows that a similar non-collapsing bound for finite times holds for more general flows of the form \eqref{F-flow} with $F$ concave (for interior non-collapsing) or convex (for exterior non-collapsing) provided ${\rm tr}\dot F$ is bounded.

\section{Proof of the main theorems}

In this section, we prove the main theorems using the maximum principle. We first derive a differential inequality of $\bar{k}(x,t)$ (in the viscosity sense). As we said in the introduction, $\bar{k}(x,t)\geq \kappa_{\max}(x,t)$, and either there exists $\bar{y}\in M\setminus\{x\}$ such that $\bar{k}(x,t)=k(x,\bar{y},t)$, or there exists a unit vector $v_0\in T_xM$ such that $\bar{k}(x,t)=h_{(x,t)}(v_0,v_0)=\kappa_{\max}(x,t)$. In the following we will treat the two cases separately.

In the case where $\bar{k}(x,t)=k(x,\bar{y},t)$ for some $\bar{y}\in M\setminus\{x\}$, we choose local normal coordinates around $x$ and $\bar{y}$. To simplify notation we denote $\omega=\f 1d(X(y,t)-X(x,t))$ and write $\pt_i^x=\f{\pt X}{\pt x^i}(x), \pt_i^y=\f{\pt X}{\pt y^i}(y)$. Then at $(x,\bar{y},t)$,
\begin{align*}
    0=\f{\pt k}{\pt y^i}=-\f 2{d^2}\langle \pt_i^y,\nu(x)+kd\omega\rangle.
\end{align*}
Noting that $\langle X(\bar{y}),\nu(x)+kd\omega\rangle=0$ and $\|\nu(x)+kd\omega\|^2=1$, we have
\begin{align}\label{nu-y}
    \nu(\bar{y})=\nu(x)+kd\omega.
\end{align}
On the other hand, a straightforward calculation gives $\langle X(\bar{y}),\pt_i^x-2\langle\pt_i^x,\omega\rangle\omega\rangle=0$ and
$\langle\nu(\bar{y}),\pt_i^x-2\langle\pt_i^x,\omega\rangle\omega\rangle=0.$ So we conclude that the plane spanned by $\pt_i^x-2\langle\pt_i^x,\omega\rangle\omega,i=1,\cdots,n$ coincides with the plane spanned by $\pt_i^y, i=1,\cdots,n$. By a suitable choice of the coordinates system near $\bar{y}$, we can arrange that
\begin{align}\label{coord-y}
    \pt_i^y=\pt_i^x-2\langle\pt_i^x,\omega\rangle\omega,\quad i=1,\cdots,n.
\end{align}

We first calculate the first spatial derivatives of $k$ at $(x,\bar{y},t)$.
\begin{align}
    (\f{\pt}{\pt x^i}+\f{\pt}{\pt y^i})k=&\f 2{d^2}\left(\langle \pt_i^x-\pt_i^y,\nu(x)+kd\omega\rangle-d\langle \omega,h_i^{p}(x)\pt_p^x\rangle\right)\nonumber\\
    =&\f 2{d^2}\left(\langle \pt_i^x-\pt_i^y,\nu(\bar{y})\rangle-d\langle \omega,h_i^{p}(x)\pt_p^x\rangle\right)\nonumber\\
    =&\f 2d(k-\kappa_i)\langle \pt_i^x,\omega\rangle,\label{1st-diff}
\end{align}
where we $\kappa_i$ denotes the principal curvatures at $(x,t)$.

By the homogeneity of $F$, we have $F(x)=\dot{F}^{ij}(x)h_{ij}(x)$, here $\dot{F}^{ij}$ is the derivative of $F$ with respect to the components $h_{ij}$ of the second fundamental form. We assume $F$ is concave, then for any $y\neq x$ we have $\dot{F}^{ij}(x)h_{ij}(y)\geq F(y)$ (see \cite[Lemma 5]{ALM}). Since the proof is easy, we include it here for convenience: By the concavity of $F$, we have
\begin{align*}
    F(y)\leq& F(x)+\dot{F}^{ij}(x)(h_{ij}(y)-h_{ij}(x))\\
    =&F(x)+\dot{F}^{ij}(x)h_{ij}(y)-F(x)\\
    \leq&\dot{F}^{ij}(x)h_{ij}(y),
\end{align*}
as claimed, where the equality used the Euler relation $F(x)=\dot{F}^{ij}(x)h_{ij}(x)$  by the homogeneity of $F$. The inequality is reversed when $F$ is convex.

Then we compute the second spatial derivatives of $k$ at $(x,\bar{y},t)$.
\begin{align*}
    &\dot{F}^{ij}(x)(\f{\pt}{\pt x^i}+\f{\pt}{\pt y^i})(\f{\pt}{\pt x^j}+\f{\pt}{\pt y^j})\biggr|_{(x,\bar{y},t)}k\\
    =&\f 2{d^2}\dot{F}^{ij}(x)\biggl(\langle h_{ij}(\bar{y})\nu(\bar{y})-h_{ij}(x)\nu(x)+c\delta_{ij}d\omega,\nu(\bar{y})\rangle\\
    &\quad +\langle \pt_i^x-\pt_i^y,2h_j^p(x)\pt_p^x+2d\omega(\f{\pt}{\pt x^j}+\f{\pt}{\pt y^j})k+k(\pt_j^y-\pt_j^x)\rangle\\
    &\quad -d\langle\omega,\nabla_jh_i^p(x)\pt_p^x\rangle+d\langle\omega,ch_{ij}(x)X(x)+h_i^p(x)h_{pj}(x)\nu(x)\rangle\biggr)\\
    \geq& \f 2{d^2}\biggl(\langle F(\bar{y})\nu(\bar{y})-F(x)\nu(x),\nu(\bar{y})\rangle+\f{d^2}2{\rm tr}(\dot{F})ck-d\langle \omega,\nabla F(x)\rangle\\
    &\quad +4(k-\kappa_i)\langle\pt_i,\omega\rangle\langle\pt_j,\omega\rangle\dot{F}^{ij}(x)\biggr)-cF(x)-\dot{F}^{ij}(x)h_i^p(x)h_{pj}(x)k,
\end{align*}
where we used \eqref{defn-k},\eqref{nu-y},\eqref{coord-y}, \eqref{1st-diff} and the inequality $\dot{F}^{ij}(x)h_{ij}(\bar{y})\geq F(\bar{y})$. Here ${\rm tr}(\dot{F})$ denotes the trace of the matrix $\dot{F}^{ij}$.

Noting that the evolution of $\nu(x,t)$ is given by
\begin{align*}
    \f{\pt\nu}{\pt t}(x,t)=\nabla F(x,t)+cF(x,t)X(x,t),
\end{align*}
the time derivative of $k$ at $(x,\bar{y},t)$ can be calculated as
\begin{align*}
    \f{\pt k}{\pt t}=&\f 2{d^2}\biggl(\langle F(\bar{y})\nu(\bar{y})-F(x)\nu(x),\nu(\bar{y})\rangle-d\langle \omega,\nabla F(x)+cF(x)X(x)\rangle\biggr)
\end{align*}

Since $\bar{k}$ is in general not smooth, we prove that $\bar{k}$ satisfies a differential inequality in a viscosity sense: For an arbitrary $C^2$ function $\phi$ which touches $\bar{k}$ from above on a neighbourhood of $(x,t)$ in $M\times [0,t]$, with equality at $(x,t)$, we prove the differential inequality for $\phi$ at $(x,t)$. From $\phi(x,t)=\bar{k}(x,t)=k(x,\bar{y},t)$, and $\phi(x',t')\geq k(x',y',t')$ for all points $x'$ near $x$, $y'\neq \bar{y}$ and earlier time $t'\leq t$, we conclude that at $(x,t)$
\begin{align*}
    \biggl(\f{\pt}{\pt t}-\dot{F}^{ij}\nabla_i\nabla_j\biggr)\phi(x,t)\leq&\biggl(\f{\pt}{\pt t}-\dot{F}^{ij}(\f{\pt}{\pt x^i}+\f{\pt}{\pt y^i})(\f{\pt}{\pt x^j}+\f{\pt}{\pt y^j})\biggr)\biggr|_{(x,\bar{y},t)}k\\
    \leq &\left(\dot{F}^{ij}h_i^ph_{pj}-tr(\dot{F})c\right)k(x,\bar{y},t)+2cF(x,t)\\
    &\quad -\f 8{d^2}(k(x,\bar{y},t)-\kappa_i(x,t))\langle\pt_i,\omega\rangle\langle\pt_j,\omega\rangle\dot{F}^{ij}(x,t)\\
    \leq&\left(\dot{F}^{ij}h_i^ph_{pj}-tr(\dot{F})c\right)\phi(x,t)+2cF(x,t),
\end{align*}
where we used $k(x,\bar{y},t)=\bar{k}(x,t)\geq \kappa_i(x,t)$ and the fact that the matrix $\dot{F}^{ij}$ is positive definite.

We now consider the case $\bar{k}(x,t)=h_{(x,t)}(v_0,v_0)=\kappa_{\max}(x,t)$. We define a smooth unit vector field $v$ near $(x,t)$ by choosing $v(x,t)=v_0$, extending in space by parallel translation along geodesics, and extending in time by solving $\f{\pt v}{\pt t}=F\mathcal{W}(v)$, where $\mathcal{W}$ is the Weingarten map. We need the following lemma about the evolution equation for the second fundamental form.
\begin{lem}[\cite{An-94}]\label{lem-3}
Under the curvature flow \eqref{F-flow} in $N^{n+1}(c)$, we have
\begin{align*}
    \f{\pt h_i^j}{\pt t}=&\dot{F}^{kl}\nabla_k\nabla_l h_i^j+\ddot{F}^{kl,pq}\nabla_ih_{kl}\nabla^jh_{pq}+(\dot{F}^{kl}h_k^ph_{pl}-tr(\dot{F})c)h_i^j+2cF\delta_i^j.
\end{align*}
\end{lem}
By the concavity of $F$, the second term on the righthand side is non-positive. Since $\phi=\bar{k}=h(v,v)$ at the point $(x,t)$, and $\phi\geq\bar{k}\geq h(v,v)$ at nearby points and earlier times, we have
\begin{align*}
    \biggl(\f{\pt}{\pt t}-\dot{F}^{ij}\nabla_i\nabla_j\biggr)\phi(x,t)\leq& \biggl(\f{\pt}{\pt t}-\dot{F}^{ij}\nabla_i\nabla_j\biggr)\biggr|_{(x,t)}h(v,v)\\
    \leq&\left(\dot{F}^{ij}h_i^ph_{pj}-tr(\dot{F})c\right)\phi(x,t)+2cF(x,t).
\end{align*}

So we conclude that the function $\bar{k}(x,t)$ satisfies the following differential inequality in a viscosity sense:
\begin{align}\label{evl-k}
    \biggl(\f{\pt}{\pt t}-\dot{F}^{ij}\nabla_i\nabla_j\biggr)\bar{k}\leq&\left(\dot{F}^{ij}h_i^ph_{pj}-tr(\dot{F})c\right)\bar{k}+2cF.
\end{align}

We now complete the proof of Theorem \ref{sphere}. Recall that under the curvature flow \eqref{F-flow}, the evolution of the speed $F$ is given by (see \cite{An-94})
\begin{align}
    \left(\f{\pt}{\pt t}-\dot{F}^{ij}\nabla_i\nabla_j\right)F=&\left(\dot{F}^{ij}h_i^ph_{pj}+tr(\dot{F})c\right)F.
\end{align}
When $F$ is positive, we define $\varphi(t)=e^{2nt}(\sup_{x\in M}\frac{\bar{k}}{F}-\frac 1n)$ for each time $t$. We show that $\varphi(t)$ is non-increasing in $t$. It suffices to prove that $\bar{k}(x,t)-(\f 1n+e^{-2nt}\varphi(t_0)+\epsilon e^{t-t_0})F(x,t)\leq 0$ for any $t_0\in [0,T)$, $t\in [t_0,T)$ and any $\epsilon>0$. Taking $\epsilon\ra 0$ then gives $\bar{k}(x,t)-(\f 1n+e^{-2nt}\varphi(t_0))F(x,t)\leq 0$ and therefore $\varphi(t)\leq \varphi(t_0)$ for $t_0\leq t$.

At time $t_0$, we have $\bar{k}(x,t_0)-(\f 1n+e^{-2nt_0}\varphi(t_0)+\epsilon)F(x,t_0)\leq -\epsilon F(x,t_0)< 0$ for all $x$. So if $\bar{k}-(\f 1n+e^{-2nt}\varphi(t_0)+\epsilon e^{t-t_0})F$ does not remain negative for $t>t_0$, there exists a first time $t_1>t_0$ and some point $x_1\in M$ such that $\bar{k}-(\f 1n+e^{-2nt}\varphi(t_0)+\epsilon e^{t-t_0})F$ is non-positive on $M\times [t_0,t_1]$ but $\bar{k}(x_1,t_1)-(\f 1n+e^{-2nt_1}\varphi(t_0)+\epsilon e^{t_1-t_0})F(x_1,t_1)=0$, i.e., the function $\phi(x,t)=(\f 1n+e^{-2nt}\varphi(t_0)+\epsilon e^{t-t_0})F(x,t)$ touches $\bar{k}(x,t)$ from above in $M\times [t_0,t_1]$, with equality at $(x_1,t_1)$. Since $\bar{k}(x,t)$ satisfies the differential inequality \eqref{evl-k} in a viscosity sense, we have that at the point $(x_1,t_1)$ (note that in sphere case, $c=1$.)
\begin{align*}
    0\leq &-\left(\f{\pt}{\pt t}-\dot{F}^{ij}\nabla_i\nabla_j\right)\phi+(\dot{F}^{ij}h_i^ph_{pj}-tr(\dot{F}))\phi+2F\\
    =&(2ne^{-2nt_1}\varphi(t_0)-\epsilon e^{t_1-t_0})F+2F\\
    &\quad -(\f 1n+e^{-2nt_1}\varphi(t_0)+\epsilon e^{t_1-t_0})(\dot{F}^{ij}h_i^ph_{pj}+tr(\dot{F}))F\\
    &\quad +(\f 1n+e^{-2nt_1}\varphi(t_0)+\epsilon e^{t_1-t_0})(\dot{F}^{ij}h_i^ph_{pj}-tr(\dot{F}))F\\
    \leq&(2ne^{-2nt_1}\varphi(t_0)-\epsilon e^{t_1-t_0})F+2F\\
    &\quad -2n(\f 1n+e^{-2nt_1}\varphi(t_0)+\epsilon e^{t_1-t_0})F\\
    =&-(2n+1)\epsilon e^{t_1-t_0}F<0,
\end{align*}
where the second inequality used $tr{\dot{F}}\geq n$, which is due to the concavity of $F$ and $F(1,\cdots,1)=n$. This contradiction implies that $\bar{k}(x,t)-(\f 1n+e^{-2nt}\varphi(t_0)+\epsilon e^{t-t_0})F(x,t)$ remains negative. In the case where $F$ is convex and positive, we consider $\underline{k}$ instead of $\bar{k}$, all the inequalities are reversed. Then we can apply a similar argument to complete the proof of Theorem \ref{sphere}.

The case 1 of Theorem \ref{hbolic} follows from a similar argument by setting $F(x,t)=H(x,t)$ and $c=-1$. To show the second case in Theorem \ref{hbolic}, we note that under the mean curvature flow in $\H^{n+1}$, the condition $H(x,t)>n$ is preserved \cite{Hu86}. We need to show that  $\psi(t)=\sup_{x\in M}\f{\bar{k}}{H-n}$ is non-increasing in $t$. As in the proof of Theorem \ref{sphere}, it suffices to show that $\bar{k}(x,t)-(H(x,t)-n)(\psi(t_0)+\epsilon e^{t-t_0})\leq 0$ for any $t_0\in [0,T)$, $t\in [t_0,T)$ and any $\epsilon>0$. The remaining argument is similar.

%=============================================================

\bibliographystyle{Plain}

\begin{thebibliography}{10}

\bibitem{An-94} Ben Andrews, {\it Contraction of convex hypersurfaces in Riemannian spaces}, J. Diff. Geom., \textbf{39}(1994), 407-431.

\bibitem{An} \underline{\qquad\quad}, {\it Non-collapsing in mean-convex mean curvature flow}, arXiv:1108.0247v1.

\bibitem{AB1} Ben Andrews and Paul Bryan, {\it Curvature bounds by isoperimetric comparison for normalized Ricci flow on the two-sphere,}
Calc. Var. Partial Differential Equations \textbf{39} (2010), no. 3-4, 419-428.

\bibitem{AB2} \underline{\qquad\quad}, {\it Curvature bound for curve shortening flow via distance comparison
and a direct proof of Grayson's theorem,} J. Reine Angew. Math. \textbf{653}(2011), 179-187.

\bibitem{AB3} \underline{\qquad\quad}, {\it A comparison theorem for the isoperimetric profile under curve shortening flow,} Comm.
Analysis and Geometry, \textbf{19} (2011), 503-530.

\bibitem{ALM} Ben Andrews, Mat Langford, and James McCoy, {\it Non-collapsing in fully nonlinear curvature flow}, arXiv: 1109.2200v1.

\bibitem{AL} Ben Andrews, Haizhong Li, {\it Embedded constant mean curvature tori in the three-sphere}, arXiv:1204.5007v2.

\bibitem{Bren} Simon Brendle, {\it Embedded minimal tori in $\SS^3$ and the Lawson conjecture}, arXiv: 1203.6597v2.

\bibitem{Ha1} Richard S. Hamilton, {\it Harnack estimate for the mean curvature flow}, J. Differential Geom. \textbf{41}(1995), no.1, 215-226.

\bibitem{Ha2} \underline{\qquad\quad} , {\it An isoperimetric estimate for the Ricci flow on the two-sphere}, Modern methods in complex
analysis (Princeton, NJ, 1992), Ann. of Math. Stud., vol. 137, Princeton Univ. Press, Princeton, NJ, 1995, pp. 191-200.

\bibitem{Ha3}\underline{\qquad\quad}, {\it Isoperimetric estimates for the curve shrinking flow in the plane}, Modern methods in complex analysis (Princeton, NJ, 1992), Ann. of Math. Stud., vol. 137, Princeton Univ. Press, Princeton, NJ, 1995, pp. 201-222.

\bibitem{Hu86} Gerhard Huisken, {\it Contracting convex hypersurfaces in Riemannian manifolds by their mean curvature}, Invent. math. \textbf{84}(1986), 463-480.

\bibitem{Hu87} \underline{\qquad\quad}, {\it Deforming hypersurfaces of the sphere by their mean curvature}, Math. Z. \textbf{195}(1987), 205-219.

\bibitem{Hu98}\underline{\qquad\quad}, {\it A distance comparison principle for evolving curves,} Asian J. Math. \textbf{2}(1998), no.1, 127-133.

\bibitem{SW} Weimin Sheng and Xu-Jia Wang, {\it Singularity profile in the mean curvature flow}, Methods Appl. Anal. \textbf{16}(2009), no.2, 139-155.

\bibitem{W} Brian White, {\it The size of the singular set in mean curvature flow of mean-convex sets}, J. Amer. Math. Soc. \textbf{13}(2000), no.3, 665-695.

\end{thebibliography}

\end{document}